\documentstyle{amsppt}
\tolerance 3000
\pagewidth{5.5in}
\vsize7.0in
\magnification=\magstep1
\widestnumber \key{AAAAAAA}
\topmatter
\author Alex Iosevich    
\endauthor
\title Fourier transform, $L^2$ restriction theorem, and scaling    
\endtitle 
\abstract We show, using a Knapp-type homogeneity argument, that the
$(L^p, L^2)$ restriction theorem implies a growth condition on the
hypersurface in question. We further use this result to show that the
optimal $(L^p, L^2)$ restriction theorem implies the sharp isotropic
decay rate for the Fourier transform of the Lebesgue measure carried by
compact convex finite hypersurfaces. 
\endabstract 
\endtopmatter 
\document 

\head Section 0: Introduction \endhead 
\vskip.125in  

Let $S$ be a smooth compact finite type hypersurface. Let 
$F_S(\xi)=\int_S e^{-i\langle x,\xi \rangle} d\sigma(x)$, where $d\sigma$ 
denotes the Lebesgue measure on $S$. Let ${\Cal R}f=\hat{f}{|}_S$, where 
$\hat{f}(\xi)$ denotes the standard Fourier transform of $f$. It is well 
known (see \cite{T}, \cite{Gr}) that if $|F_S(\xi)| \leq C{(1+|\xi|)}^{-r}$, 
$r>0$, then ${\Cal R}: L^p({\Bbb R}^n) \rightarrow L^2(S)$ for $p \leq 
\frac{2(r+1)}{r+2}$. A natural question to ask is, does the boundedness of  
${\Cal R}: L^{\frac{2(r+1)}{r+2}}({\Bbb R}^n)\rightarrow L^2(S)$, $r>0$, 
imply that $|F_S(\xi)| \leq C{(1+|\xi|)}^{-r}$? In this paper we will show
that this is indeed the case if $S$ is a smooth convex finite type hypersurface
in the sense that the order of contact with every tangent line is finite. (See
\cite{BNW}). 

Under a more general condition, called the finite polyhedral type assumption,  
(see Definition 1 below), we will show that the $(L^p, L^2)$ restriction  
theorem with $p=\frac{2(r+1)}{r+2}$ implies that $|B(x,\delta)| \leq C{\delta} 
^r$, where $B(x, \delta)=\{y: dist(y, T_x(S)) \leq \delta\}$, where $T_x(S)$ 
denotes the tangent hyperplane to $S$ at $x$.  

Our plan is as follows. We will first use a variant of the Knapp 
homogeneity argument to show that if $S$ satisfies the finite polyhedral 
type condition and ${\Cal R}: L^p({\Bbb R}^n) \rightarrow L^2(S)$, with 
$p=\frac{2(r+1)}{r+2}$, then $|B(x,\delta)|\leq C{\delta}^r$ for each $x$. 
If the surface is, in addition, convex and finite type, then the result due to
Bruna, Nagel, and Wainger (see \cite{BNW}) implies that $|F_S(\xi)| \leq 
C{(1+|\xi|)}^{-r}$. If the surface is not convex finite type, then we do 
not, in general, know how to conclude that $|B(x, \delta)| \leq C{\delta}^r$
implies that $|F_S(\xi)| \leq C{(1+|\xi|)}^{-r}$. A gap remains. 

\head Section 1: Statement of Results \endhead 
\vskip.125in 

\definition{Definition 1} Let $S$ be a smooth compact hypersurface in 
${\Bbb R}^n$. Let $B^{\pi}(x,\delta)$ denote the projection of $B(x,\delta)$
onto $T_x(S)$. We say that $S$ is of finite polyhedral type if there exists
a family of polyhedra $P(x,\delta)$ such that $B^{\pi}(x,\delta) \subset 
P(x,\delta)$, $C_1|B^{\pi}(x,\delta)| \leq |P(x,\delta)| \leq 
C_2|B^{\pi}(x,\delta)|$, where $C_1, C_2$ do not depend on $\delta$, and 
the number of vertices of $P(x, \delta)$ is bounded above independent of 
$\delta$, where $\chi_{P(x,\delta)}$ denotes the characteristic  
function of $P(x, \delta)$. 
\enddefinition  

\remark{Remark} The motivation for the definition of finite polyhedral 
type is the standard homogeneity argument due to Knapp. In order to prove 
the sharpness of the $(L^p, L^2)$ restriction theorem for hypersurfaces 
with non-vanishing Gaussian curvature, Knapp approximated such a surface
with a box with side-lengths $(\delta, \dots, \delta, {\delta}^2)$, 
$\delta$ small. He then took $f_{\delta}$ to be the inverse Fourier 
transform of the characteristic function of that box. It is not hard to see
that ${||{\Cal R}f||}_2 \approx {\delta}^{\frac{n-1}{2}}$. On the other hand,
using the fact that the Fourier transform of the box in question is $\approx
\frac{sin({\delta}^2x_n)}{x_n} \Pi_{j=1}^{n-1} \frac{sin(\delta x_j)}{x_j}$, 
it is not hard to see that ${||f_{\delta}||}_p \approx 
{\delta}^{\frac{n+1}{p'}}$. It follows that $p \leq \frac{2(n+1)}{n+3}$, 
which is the known positive result due to Stein and Tomas. The crucial 
part of this calculation is the approximation of the surface with a box 
with appropriate dimensions. Definition 1 and Theorem 2 are generalizations 
of this phenomenon.  

It should also be noted that it is not hard to see that the hypersurface 
$S=\{x: x_3=x_1x_2\}$ does not satisfy the finite polyhedral type condition. 
Thus, it makes sense to think of the finite polyhedral type condition as a
generalization of convexity. 
\endremark 

\proclaim{Theorem 2} Let $S$ be of finite polyhedral type. Consider the 
estimates 
$$ |F_S(\xi)| \leq C{(1+|\xi|)}^{-r}, \tag*$$ 
$$ {\Cal R}: L^p({\Bbb R}^n) \rightarrow L^2(S), \ \ \ 
p\leq \frac{2(r+1)}{r+2}, \tag**$$ and 
$$ |B(x,\delta)| \leq C{\delta}^r, \tag***$$ for each $x$.   

Then $(*)$ implies $(**)$ and $(**)$ implies $(***)$. Further, $(***)$ implies $(*)$  
if $S$ is in addition convex and finite type. 
\endproclaim 

The fact that $(*)$ implies $(**)$ is essentially the Stein-Tomas restriction
theorem. (See \cite{T}, \cite{Gr}). 
The fact that $(***)$ implies $(*)$ in the case of convex finite
type hypersurfaces is due to Bruna, Nagel, and Wainger. (See \cite{BNW}). So 
it remains to prove that $(**)$ implies $(***)$, and that convex finite 
type hypersurfaces are of finite polyhedral type. (See Theorems 3 and 4 below). 
 
\proclaim{Theorem 3} Let $S=\{x: x_n=Q(x')+R(x')+c\}$, where 
$x'=(x_1, \dots, x_{n-1})$, $Q\in C^{\infty}$ is mixed homogeneous in the  
sense that there exist integers $(a_1, \dots, a_{n-1})$, $a_j\ge 1$, such that 
$Q(t^{\frac{1}{a_1}}x_1, \dots, t^{\frac{1}{a_{n-1}}}x_{n-1})=tQ(x')$, 
$Q(x')\not=0$ if $x'\not=(0, \dots, 0)$, $R\in C^{\infty}$ is the remainder 
in the sense that $\lim_{t \rightarrow 0} 
\frac{R(t^{\frac{1}{a_1}}x_1, \dots, t^{\frac{1}{a_{n-1}}}x_{n-1})}{t}=0,$ and
$c$ is a constant. Then $S$ is of finite polyhedral type.  
\endproclaim   

Theorem 3 implies that convex finite type hypersurfaces are of finite 
polyhedral type via the following representation result due to Schulz. (See 
\cite{Sch}. See also \cite{IS2}). 

\proclaim{Theorem 4} Let $\Phi \in C^{\infty}({\Bbb R}^{n-1})$ be a convex
finite type function such that $\Phi(0, \dots,0)=0$ and 
$\nabla \Phi(0, \dots, 0)=(0, \dots, 0)$. Then, after perhaps applying a 
rotation, we can write $\Phi(y)=Q(y)+R(y)$, where $Q$ is a convex polynomial,
mixed homogeneous in the sense of Theorem 3, and $R$ is the remainder 
in the sense of Theorem 3. 
\endproclaim 
\vskip.25in 

\head Section 2: $(**) \rightarrow (***)$ \endhead 
\vskip.125in 

Locally, $S$ is a graph of a smooth function $\Phi$, 
such that $\Phi(0, \dots, 0)=0$, and $\nabla \Phi(0, \dots, 0)=
(0, \dots, 0)$. If we consider a sufficiently small piece of our hypersurface,
$B^{\pi}(0, \dots, 0, \delta)=\{y \in K: \Phi(y)\leq \delta\}$, where $K$ is a 
compact set in ${\Bbb R}^{n-1}$ containing the origin, and, without 
loss of generality, $\Phi(y) \ge 0$. Since $|B^{\pi}(x,\delta)| \approx 
|B(x,\delta)|$, it suffices to show that 
$|\{y \in K: \Phi(y) \leq \delta \}| \leq C{\delta}^{r}$. 

Let $f_{\delta}$ be a function such that $\hat{f_{\delta}}$ is the 
characteristic function of the set $\{(x',x_n): x' \in P_{\delta} \ \ 
; \ 0\leq x_n \leq \delta\}$, where $P_\delta$ is the polyhedron containing
the set $\{x': \Phi(x') \leq \delta\}$ given by the definition of finite 
polyhedral type. 

Let's assume for a moment that ${||f_{\delta}||}_p \leq C 
{(\delta |P_{\delta}|)}^{1/p'}$. Since the 
restriction theorem holds, we must have ${||{\Cal R}f_{\delta}||}_2 \leq 
C{||f_{\delta}||}_p,$ which implies that $|P_{\delta}| \leq C{\delta}^
{\frac{2(p-1)}{2-p}}$. Since $p=\frac{2(r+1)}{r+2}$, it follows that 
$|P_{\delta}| \leq C{\delta}^r$. By the  definition of $P_{\delta}$ it follows 
that $|B(0, \dots, 0, \delta)|=|B^{\pi}(0, \dots, 0, \delta)|=
|\{y \in k: \Phi(y) \leq \delta \}| \leq 
C{\delta}^r$. 

This completes the proof provided that we can show that 
${||f_{\delta}||}_p \leq C{(\delta |P_{\delta}|)}^{1/p'}$. More generally, 
we will show that if $P$ is a polyhedron in ${\Bbb R}^n$, then 
${||\widehat{\chi}_{P}||}_p \leq C{|P|}^{1/p'}$, where $C$ depends on the
dimension and the number of vertices of $P$ 
and $|P|$ denotes the volume of $P$. We give the argument in 
two dimensions, the argument in higher dimensions being similar. Break up
$P$ as a union of disjoint (up to the boundary) triangles $t_j$, $j=1, \dots
,N$. Since $\chi_{P}(x)=\sum_{j} \chi_{t_j}(x)$, it suffices to carry out
the argument for $\chi_P$, where $P$ is assumed to be a triangle. Since 
translations don't contribute anything in this context, we may assume that
one of the vertices of the triangle is at the origin. Break up this
triangle, if necessary, into two right triangles. Refine the original
decomposition so that it consists of right triangles. Rotate the right
triangle so that it is in the first quadrant and one of the sides is on
the $x_1$-axis. We now apply a linear transformation mapping this
triangle (denoted by $P'$) into the triangle with the endpoints 
$(0,0)$, $(1,0)$ and
$(1,1)$. It is easy to check by an explicit computation that the Fourier
transform of the characteristic function of this triangle has the $L^p$
norm (crudely) bounded by $2^{\frac{1}{p}}$. 

Let $T$ denote the linear 
transformation taking the triangle $P'$ to the unit triangle above. 
We see that
$$ \widehat{\chi}_{TP'}(\xi)=|T| \widehat{\chi}_{P'}(T^t\xi),$$ so 
$${||\widehat{\chi}_{TP'}||}_p={|T|}^{1/p'} {||\widehat{\chi}_{P'}||}_p.$$ 

Since $|T|=\frac{1}{|P'|}$, we see that ${||\widehat{\chi}_{t_j}||}_p \leq 
C{|t_j|}^{1/p'}$, where the $t_j$'s are the triangles from the (refined)
original decomposition. Adding up the estimates we get 
$${||\widehat{\chi}_{P}||}_p \leq C \sum_{j=0}^{N} {|t_j|}^{1/p'} \leq 
CN {(\sum_{j=0}^N |t_j|)}^{1/p'}=CN{|P|}^{1/p'}.$$  

In higher dimensions the proof is virtually identical with triangles 
replaced by $n-1$ dimensional simplices, i.e the convex hull of $n$ points 
that are not contained in any $(n-2)$ dimensional plane. 

Since we have assumed that the number of vertices of $P_{\delta}$ is bounded 
above, it follows that ${||f_{\delta}||}_p \leq  
C{(\delta |P_{\delta}|)}^{1/p'}$, as desired. 
\vskip.25in 

\head Section 3: Proof of Theorem 3 \endhead 
\vskip.125in 

As before, it is enough to consider the set $B^{\pi}(0, \dots, 0)=
\{y: Q(y)+R(y) \leq \delta\}$. It will be clear from the proof below that if 
we shrink the support sufficiently, then $B^{\pi}(0, \dots, 0) \approx  
B_{Q}^{\delta}=\{y: Q(y) \leq \delta\},$ due to our assumptions on the 
remainder term $R$.  

Let $\frac{n-1}{m}=\frac{1}{a_1}+
\dots + \frac{1}{a_{n-1}}$. 
Our plan is as follows. We first prove that $|B_{Q}^{\delta}| \approx 
{\delta}^{\frac{n-1}{m}}$. Then, we will find a polyhedron of suitable area that
contains the set $B_{Q}^{1}$. We shall obtain the polyhedra for all values of
$\delta$ by homogeneity.  

Going into polar coordinates,
$x_1=s^{\frac{m}{a_1}}\omega_1, \dots, x_{n-1}=s^{\frac{m}{a_{n-1}}}$, 
$\omega=(\omega_1, \dots, \omega_{n-1}) \in S^{n-2}$, we see that
$\int_{B_{Q}^{\delta}} dy=\int_{S^{n-2}} 
\int_{0}^{s^{\frac{1}{m}}Q^{-\frac{1}{m}}(\omega)} s^{n-2}ds d\omega =
{\delta}^{\frac{n-1}{m}} \int_{S^{n-2}} Q^{-\frac{n-1}{m}}(\omega) d\omega=
C_{Q}{\delta}^{\frac{n-1}{m}}$. This proves that $|B_Q^{\delta}| \approx 
{\delta}^{\frac{n-1}{m}}$. 

We now find a box $P_1$ with sides parallel to the coordinate axes, 
such that $B_{Q}^{1} \subset P_1$, and $|P_1|=cC_Q$, where $c>1$. 
Let $Q_P$ be a mixed homogeneous function of degree 
$(a_1, \dots, a_{n-1})$ defined by the condition $\{y:Q_P(y)=1\}=\partial P_1$, 
where $\partial P_1$ denotes the boundary of $P_1$. Let $P_{\delta}$ be the  
polyhedron such that the boundary $\partial P_{\delta}=\{y: Q_P(y)=\delta\}$.
It is not hard to see that $B_{Q}^{\delta} \subset P_{\delta}$. Moreover, 
$|P_{\delta}|=cC_Q{\delta}^{\frac{n-1}{m}}\approx |B_Q^{\delta}|$ by the 
calculation made in the previous paragraph. This completes the proof of 
Theorem 3. 
\newpage 

\head References \endhead 
\vskip.125in 

\ref \key BNW \by J. Bruna, A. Nagel, and S. Wainger  
\paper Convex hypersurfaces and Fourier transform  
\jour Annals of Math. \pages 333-365 \vol 127 \yr 1988 \endref 

\ref \key Gr \by A. Greenleaf \pages 519-537 
\paper Principal curvature and harmonic analysis 
\yr 1982 \vol30 \jour Indiana Math J. \endref 

\ref \key IS1 \by A. Iosevich and E. Sawyer 
\paper Oscillatory integrals and maximal averaging operators associated to 
homogeneous hypersurfaces \yr 1996 \jour Duke Math J. \vol 82 \endref 

\ref \key IS2 \by A. Iosevich and E.Sawyer 
\paper Maximal averages over surfaces  
\jour (to appear in Adv. in Math.) \yr 1996 \endref 

\ref \key Sch \by H. Schulz 
\paper Convex hypersurfaces of finite type and the asymptotics of their
Fourier transforms \vol 40 \yr 1991 \jour Indiana Univ. Math. J. \endref

\ref \key St \by E. M. Stein \paper Harmonic Analysis in ${\Bbb R}^n$ 
(Proc. Conf. De Paul University) \jour Math. Assoc. Am. Studies in Math
\yr 1976 \vol 13 \endref 

\ref \key T \by P. Tomas \paper A restriction theorem for the Fourier 
transform \jour Bull. AMS \vol 81 \yr 1975 \pages 477-478 \endref 

\enddocument